\documentclass[a4paper,10pt]{article}
\usepackage{ifthen}
\usepackage{makeidx}
\usepackage{amsmath,amsfonts,amssymb}
\usepackage[nohug 
,fixed]{diagrams}
\usepackage{epsfig}

\newcommand{\pleinepage}%
{\setlength{\oddsidemargin}{0in}\setlength{\textwidth}{6.26in}\setlength{\topmargin}{0in}\setlength{\textheight}{8.7in}}
\newcommand{\grossepage}%
{\setlength{\oddsidemargin}{-0.5cm}\setlength{\textwidth}{17.5cm}\setlength{\topmargin}{-2cm}\setlength{\textheight}{26cm}}

\newcommand{\ie}{\emph{i.e.\ }}

\newcommand{\cf}{cf.\ }

\newcommand{\romaindroit}[1]{\textrm{\textup{#1}}}

\renewcommand{\epsilon}{\varepsilon}
\renewcommand{\phi}{\varphi}

\newcommand{\ZZ}{\mbox{Z\hspace{-0.3em}Z}}

\newcommand{\RR}{\mbox{I\hspace{-0.15em}R}}
\newcommand{\CC}{\hspace{0.15em}\mbox{l\hspace{-0.47em}C}}
\newcommand{\TT}{\mbox{T\hspace{-0.45em}T}}

\renewcommand{\ZZ}{\mathbb{Z}}

\renewcommand{\RR}{\mathbb{R}}
\renewcommand{\CC}{\mathbb{C}}
\renewcommand{\TT}{\mathbb{T}}

\newcommand{\inv}{^{-1}}
\newcommand{\inverse}[1]%
{\frac{1}{#1}}

\newcommand{\supp}{\hbox{\romaindroit{Supp}\,}}

\newcommand{\restr}[2]{\left. #1 \right|_{#2}}

\newcommand{\pensebete}[1]%
{\typeout{attention ! oubli de << #1 >> !}
\begin{center} {\bf #1} \end{center}}

\newcommand{\Continue}{\romaindroit{C}}

\newcommand{\Ci}{{\Continue^\infty}}%
\newcommand{\Cc}{{\Continue_\textrm{c}}}%
\newcommand{\Czero}{{\Continue_0}}%
\newcommand{\Cci}{{\Continue^\infty_\textup{c}}}%
\newcommand{\Cizero}{{\Continue_0^\infty}}%
\newcommand{\cstar}{{\ensuremath{\Continue^*}}}%

\newcommand{\fleche}[1]%
{\rTo^{#1}}
\newcommand{\fonction}[5]%
{\begin{diagram}
#2 & {} &\rTo^{#1} & {} & #3 \\
#4 & {} &\rMapsto & {} & #5
\end{diagram} }

\newcommand{\sfonction}[5]%
{$\begin{array}{ccc}#2 & {\buildrel #1 \over \rightarrow} & #3 \\#4 & \mapsto & #5 \\ \end{array}$ }
\newcommand{\ssfonction}[5]%
{\begin{array}{ccc}#2 & {\buildrel #1 \over \rightarrow} & #3 \\#4 & \mapsto & #5 \\ \end{array} }

\renewcommand{\{}{\left\lbrace\left.}
\renewcommand{\}}{\right.\right\rbrace}
\newcommand{\tq}{\ \right| \hspace{-1.1ex} \left|\ }
\newcommand{\bracket}[2]{\left\langle \left.  #1 \right|
    \hspace{-1.1ex} \left|\ #2 \right.\right\rangle}

\newcommand{\matricedeux}[4]%
{\left(\begin{array}{cc} #1 & #2 \\ #3 & #4 \end{array}\right)}

\newcommand{\matricetrois}[9]%
{\left(\begin{array}{ccc} #1 & #2 & #3 \\ #4 & #5 & #6 \\ #7 & #8 &
      #9 \end{array}\right)}

\newcommand{\matricetriangsup}[5]%
{\left(\begin{array}{cccc}
#1 & #3 & \cdots & #4 \\
 0 & \ddots & \ddots & \vdots \\
\vdots & \ddots & \ddots & #5\\
0& \cdots & 0 & #2 \end{array}\right)}

\newcommand{\matricediag}[2]%
{\left(\begin{array}{cccc}
#1 & 0 & \cdots & 0 \\
 0 & \ddots & \ddots & \vdots \\
\vdots & \ddots & \ddots & 0\\
0& \cdots & 0 & #2 \end{array}\right)}

\catcode`à=\active \defà{{\`a}}%
\catcode`â=\active \defâ{{\^a}}%
\catcode`ä=\active \defä{{\"a}}%
\catcode`é=\active \defé{{\'e}}%
\catcode`è=\active \defè{{\`e}}%
\catcode`ê=\active \defê{{\^e}}%
\catcode`ë=\active \defë{{\"e}}%
\catcode`ï=\active \defï{{\"\i}}%
\catcode`î=\active \defî{{\^\i}}%
\catcode`ô=\active \defô{{\^o}}%
\catcode`ö=\active \defö{{\"o}}%
\catcode`ù=\active \defù{{\`u}}%
\catcode`ü=\active \defü{{\"u}}%
\catcode`û=\active \defû{{\^u}}%
\catcode`ç=\active \defç{{\c c}}

\newcommand{\gloss}[1]%
 {\index{{#1}@{#1}}{\em #1}\relax}
\newcommand{\xgloss}[2]%
 {\index{{#1}@{#1}!{#2}@{#2}}{\em #1 \relax #2}\relax}
 \newcommand{\glossref}[2]%
 {\index{{#2}@{#1}}{\em #1}\relax}
 \newcommand{\xglossref}[4]%
 {\index{{#3}@{#1}!{#4}@{#2}}{\em #1 \relax #2 }\relax}

\newcommand{\defifont}[1]{\textsc{\textbf{#1}}}

\newcounter{theonum}
\makeatletter
\@addtoreset{theonum}{section}
\makeatother
\renewcommand{\thetheonum}{\thesection.\arabic{theonum}}
\newenvironment{defi}%
{\refstepcounter{theonum} \par \pagebreak[0]\medskip \noindent \defifont{Definition \thetheonum\
---\nopagebreak } \sffamily\renewcommand{\em}{\normalfont\itshape}}{
\\ \medskip}

\newenvironment{theo}[1][]%
 {\refstepcounter{theonum} \par \pagebreak[0] \medskip\noindent
 \defifont{Theorem \thetheonum{#1}\ ---\ }\nopagebreak \em}{\\ \smallskip}

\newenvironment{prop}%
{\refstepcounter{theonum} \par \pagebreak[0]\medskip\noindent 
\defifont{Proposition \thetheonum\ ---\ }\nopagebreak \em}{\\ \smallskip}

\newenvironment{lemme}%
{\refstepcounter{theonum}\par \pagebreak[0]\medskip\noindent
  \defifont{Lemma \thetheonum\ ---\ }\nopagebreak \em}{\\ \smallskip}
\newenvironment{lemma}%
{\begin{lemme}}{\end{lemme}}

\newenvironment{demo}%
{\par \medskip \pagebreak[3] %
 \noindent{\it Proof. ---\ \nopagebreak}%
\small\nopagebreak}{\hfill \nopagebreak {$\blacksquare$}
 \par \pagebreak[0] \smallskip}
\newenvironment{proof}%
{\begin{demo}}{\end{demo}}

\newcommand{\fcl}{\texttt{f}}
\newcommand{\isom}{\mathcal{I}}
\newcommand{\G}{G}
\newcommand{\GG}{\widetilde{G}}
\newcommand{\ff}{\widetilde{f}}
\newcommand{\dual}[1]{\widehat{#1}}
\newcommand{\zero}{^{(0)}}
\newcommand{\Gzero}{G\zero}
\newcommand{\acta}[2]{\alpha_{#1}\left(#2\right)}
\newcommand{\dbracket}[2]{\left\langle  #1 , #2 \right\rangle}
\newcommand{\vbracket}[2]{\left( \left.  #1 \right|
    #2 \right)}

\begin{document}
\title{Deformation quantization using groupoids. \\ Case of
toric manifolds}
\author{Fr\'ed\'eric CADET}
\date{\today}
\maketitle
\thispagestyle{empty}
\bigskip
\begin{abstract}
In the framework of \cstar-algebraic deformation quantization
we propose a notion of deformation groupoid which could apply to
known examples {\em e.g.} Connes' tangent groupoid of a manifold, its
generalisation by Landsman and Ramazan, Rieffel's noncommutative
torus, and even Landi's noncommutative 4-sphere. We
construct such groupoid for a wide class of $\TT^n$-spaces, that
generalizes the one given for $\CC^n$ by Bellissard and Vittot. In
particular, using the geometric properties of the moment map discovered
in the '80s by Atiyah, Delzant, Guillemin and Sternberg, it provides a 
\cstar-algebraic deformation quantization for all toric
manifolds, including the 2-sphere and all complex projective spaces. 
\end{abstract}

\bigskip 

\section*{Introduction}

{\em Quantization} is going from classical mechanics to quantum
mechanics. In classical hamiltonian picture, a physical system is
described by its phase space, a Poisson manifold $M$, which
Poisson bivector is denoted by $\Omega$, and associated Poisson bracket
of two observables $\fcl,\fcl'\in \Ci(M)$ by:
$$\{\fcl,\fcl'\}_\Omega\overset{def}{=} \bracket{d\fcl\otimes d\fcl'}{\Omega}.$$
On the other hand a quantum system is described, in Heisenberg's
picture, by replacing the classical algebra $\Ci(M)$ of observables by
an algebra of operators on a Hilbert space, in a linear and
involution-preserving way. Moreover the Poisson bracket $\{\fcl,\fcl'\}_\Omega$
has to be replaced by $\inverse{\imath\hbar}[f,f']$, the commutator of
$f$ and $f'$, the quantum analogs of $\fcl$ and $\fcl'$, divided by
$\imath=\sqrt{-1}$ and by $\hbar$, which is Planck's physical constant.
Dirac postulated in the '30s that this correspondance from
renormalized commutators to Poisson bracket is given by a limit
process when Planck's constant goes to zero.

In the '70s the notion of {\em deformation quantization} was
introduced in \cite{BFFLS} as an attempt to give a precise definition
to Dirac's principle. One of the ideas was to use a field of algebras
$(\mathcal{A}_\hbar)$ parametrized by a real number $\hbar$, such
that, for $\hbar=0$, the algebra $\mathcal{A}_0$ is the commutative
algebra of classical observables, and for $\hbar\neq 0$,
$\mathcal{A}_\hbar$ is the noncommutative algebra of quantum
observables. At the beginning of the '90s Rieffel proposed in
\cite{RIE2,Rie} a topological definition of the limit process making use of
the previously known notion of continuous field of \cstar-algebras
\cite{FD}. The following is a new definition close Rieffel's one, but a little more general so that there is a direct adaptation
to groupoid \cstar-algebras.
\begin{defi}
Let $(\mathcal{A}_\hbar)_{\hbar\in X}$ be a continuous field of
\cstar-algebras parametrized by a locally compact subset $X$ of $\RR$
containing $0$ as a limit point; we denote
$\mathcal{A}=\bigcup_{\hbar\in X}\mathcal{A}_\hbar$ the associated
topological bundle over $X$. Let $\mathcal{Q}$ be a sub-*-algebra of
the \cstar-algebra $\Czero(X,\mathcal{A})$ of continuous {\em
  sections} of $\mathcal{A}$. We say that $(
(\mathcal{A}_\hbar)_{\hbar\in X},\mathcal{Q})$ is a {\em deformation}
if:
\begin{itemize}
\item The space $\mathcal{Q}_0=\{f_0 \in \mathcal{A}_0 \tq f\in
  \mathcal{Q}\}$ is dense in $\mathcal{A}_0$.
\item There is a map $\mathcal{Q}_0\times \mathcal{Q}_0\rTo \mathcal{Q}_0$,
  denoted by:
$$(a,b) \rMapsto \{a,b\}_0,$$
such that, for every $f,f'\in \mathcal{Q}$, $\{f_0,f_0'\}_0$ is the
continuous extension at zero of the the continuous section $\hbar
  \rMapsto \inverse{\imath\hbar}[f_\hbar,f'_\hbar]$ defined on
  $X-\{0\}$ \ie:
$$\lim_{\pile{_{\hbar \to 0}\\ _{\hbar \neq 0}}}
\inverse{\imath\hbar}[f_\hbar,f'_\hbar]=\{f_0,f_0'\}_0.$$
\end{itemize}\
\end{defi}
Then $\mathcal{A}_0$ is a commutative \cstar-algebra and
$\mathcal{Q}_0$ is a sub-*-algebra, and moreover $\mathcal{Q}_0$ has a
structure of Poisson algebra, for
the bracket $\{.,.\}_0$.
\begin{defi}
Let $(M,\Omega)$ be a Poisson manifold. A {\em deformation of $M$} is
a deformation $((\mathcal{A}_\hbar)_{\hbar\in X},\mathcal{Q})$ endowed
with an isomorphism of \cstar-algebras~:
$$\Czero(M) \rTo^\isom \mathcal{A}_0,$$
such that:
\begin{itemize}
\item $\mathcal{Q}_0 \subset \isom(\Cizero(M))$ 
\item  for each $f,f'\in \mathcal{Q}$:
$\{f_0,f'_0\}_0=\isom\left(\{\isom\inv f_0,\isom\inv
  f_0'\}_\Omega\right).$
\end{itemize}\
\end{defi}
Then each linear *-preserving section $T$ of the canonical projection
of $\mathcal{Q}$ onto $\mathcal{Q}_0$ gives rise to a quantization:
$$\begin{diagram} \mathcal{Q}_0 & \rTo^{quantization} &
  \mathcal{Q}_\hbar \\
& \rdDashto(1,2)_T \ \luOnto(1,2) \ruOnto(1,2) \\
& \mathcal{Q}
\end{diagram}$$

Many usual \cstar-algebras can be described as \cstar-algebras
$\cstar(G)$ of a
groupoid $G$, using the construction given by Jean Renault in
\cite{REN} that
generalizes that of the (full and reduced) \cstar-algebra(s) of a
group or of a group action. A relevant fact is that a bundle of groupoids
is itself a groupoid. On the topological level Ramazan \cite{RAM,LR}
has established the following property:
\begin{theo}
Let $(G_\hbar)_{\hbar\in X}$ be a field of groupoids, and
  $G=\bigcup_{\hbar\in X} G_\hbar$ be the corresponding bundle which is
  supposed to be locally compact Hausdorff separable, having
  continuous Haar system, and such that the bundle map $G \rOnto^p X$
  is open. If $G$ is amenable \cite{ANA} then the field
  $\left(\cstar(G_\hbar)\right)_{\hbar\in X}$ possesses a structure of
  continuous field such that the algebra of continuous sections is:
$$\Czero\left(X,\bigcup_{\hbar\in X}
  \cstar(G_\hbar)\right)=\cstar(G).$$
\end{theo}
Then it is natural to set:
\begin{defi}
A locally compact Hausdorff separable amenable bundle of groupoids
$G$ over a locally compact subset $X$ of $\RR$, with
continuous Haar system is a {\em deformation groupoid} if there exists a
sub-*-algebra $\mathcal{Q}$ of sections of the continuous field
$\left(\cstar(G_\hbar)\right)_{\hbar\in X}$ 
such that $\left(\left(\cstar(G_\hbar)\right)_{\hbar\in
    X},\mathcal{Q}\right)$ is a deformation (resp. is a {\em
  deformation groupoid of a Poisson manifold $M$} if there exists $\mathcal{Q}$ and an
isomorphism $\isom$
such that $\left(\left(\cstar(G_\hbar)\right)_{\hbar\in
    X},\mathcal{Q},\isom \right)$) is a deformation  of
$M$).
\end{defi}
Finally a bundle of groupoids $G$ over $X$ (locally compact, Hausdorff,
separable, amenable, with continuous Haar system and open projection $G\rOnto ^p X$) is a deformation
groupoid of $(M,\Omega)$ if and only if:
\begin{enumerate}
\item There exists an isomorphism $\Czero(M) \rTo_\sim^\isom \cstar(G_0)$.
\item There exists a sub-*-algebra $\mathcal{Q}$ of $\cstar(G)$ such
  that:
\begin{enumerate}
\item $\mathcal{Q}_0$ is dense in  $\cstar(G_0)$,
\item $\isom\inv(\mathcal{Q}_0) \subset \Cizero(M)$,
\item $\displaystyle \forall f,f'\in \mathcal{Q},\ \lim_{\pile{ _{\hbar \to 0}\\ _{\hbar\neq 0}}} \inverse{\imath\hbar}[f_\hbar,f'_\hbar]=\isom\left(\{\isom\inv f_0,\isom\inv
  f_0'\}_\Omega\right).$
\end{enumerate}
\end{enumerate}
Moreover we will use:
\begin{prop}\label{deformation}
Let $M$ be a  manifold, without {\em a priori} Poisson structure. If
there exists a bundle of groupoids $G$ over $X$ (locally compact,
Hausdorff, separable, amenable, with continuous Haar system) 
which satistifies the previous conditions {\em 1,2a,2b}, and such that 
\begin{itemize}
\item[{\em 2c'}]
for every $f,f'\in
\mathcal{Q}$ the section $\hbar\rTo
\inverse{\imath\hbar}[f_\hbar,f'_\hbar]$ has a unique continuous extension at $\hbar=0$ in
$\mathcal{Q}_0$, which only depends on the values $f_0$
and $f_0'$, 
\end{itemize}
then $M$ admits a Poisson bracket such that $G$ is a
deformation groupoid of $M$.
\end{prop}

The first {\em example} of deformation groupoid was given by Connes in
\cite{CON} (\cf \cite{VAR} too): the tangent groupoid of a manifold
$N$ is a deformation groupoid of the cotangent bundle $T^*N$ endowed with its
canonical symplectic structure. This was generalized by Landsman and
Ramazan \cite{LR,RAM,LAN} to integrable Lie-Poisson manifolds (\ie
manifolds which are the dual of an integrable Lie algebroid, endowed whith the
canonical Poisson strucure described in \cite{DAZ}
. In 1990
Bellissard and Vittot \cite{BEL} constructed a deformation groupoid of $\CC^n$
with its canonical symplectic structure, which was not a Lie
groupoid. In the present paper, we generalize  this
construction to other $\TT^n$-spaces $M$ ($\TT^n$ denotes the
$n$-dimensional torus); it applies to the 2-sphere and, using
Delzant's results \cite{DEL}, more generally
to all toric manifolds, including all complex projective spaces,
endowed with their canonical 
Kähler structure. 

The strategy is based on the two remarks, first that the previous condition (1.)
implies by \cite[lemme 1.3 p.7]{REN2} that the ``classical'' groupoid $G_0$
must be a bundle of commutative groups\footnote{Remark that a bundle of
  commutatives groups is always amenable. Then replacing globally
  continuous fields of \cstar-algebras by fields only continuous at
  zero, one can define non globally amenable deformation groupoids}, and second that an example of
such isomorphism $\isom$ is the Fourier-Gelfand transform
$\Continue(\TT^n) \simeq \cstar(\ZZ^n)$. This paper is organized as
follows: 
\begin{enumerate}
\item  For a $\TT^n$-space $M$, we construct an isomorphism
  $\Continue(M) \simeq \cstar(G_0)$ for a convenient bundle of
  commutative groups $G_0$, under the assumption that the
  projection $M \rOnto M/\TT^n$ has a continuous section.
\item Using a second action on $M/\TT^n$ we construct then a deformation groupoid bundle $G$ over
  $\RR$, such that the fiber at 0 of $G$ is $G_0$. The sub-*-algebra $\mathcal{Q}$
  considered is made of restrictions of $\Cci$ functions on a Lie
  groupoid $\GG$ containing $G$ as subgroupoid. 
\item We investigate which Poisson structure corresponds to
  this deformation structure on $M$.
\item We verify that this construction can be applied to toric
  manifolds such that the Poisson structure obtained equals the
  original symplectic structure of the toric manifold.
\end{enumerate}
\section{Fourier-Gelfand isomorphisms for a $\TT^n$-space}
For $G$ a groupoid, we denote $\Gzero$ its units space, and $G
\pile{\rTo^r \\ \rTo_s} \Gzero$ its range and source maps. Recall the
definition \cite{REN}:
\begin{defi}
A locally compact, Hausdorff, separable groupoid $G$ is {\em \'etale}
when
the fibers of the range map: $$G^y=\{g\in G \tq r(g)=y \},\ \  y\in\Gzero,$$
are discrete, and the counting measure is a continuous Haar system.
\end{defi}
Let us recall some known properties:
\begin{itemize}
\item An open subgroupoid of an \'etale groupoid is
itself -- locally compact, Hausdorff, separable -- \'etale.
\item The maps $r,s$ are automatically open \cite{REN}. And the converse is true
  for bundle of
commutative groups: such topological bundle is étale if and only if the projection $G
\rTo^{r=s} \Gzero$ is open (\cf \cite{REN2}).
\end{itemize}
Let $M$ be a Hausdorff, locally compact, separable space endowed with
an action $\alpha$ of $\TT^n$ denoted by: $(s,x)\in \TT^n\times M \rMapsto
\acta{s}{x}\in M$. We construct a bundle of commutative
groups in the following way. Let us introduce some notations:
\begin{itemize}
\item $\Delta=M/\TT^n$ denotes the quotient space, which is Hausdorff, separable and
  locally compact, and $M\rOnto^J \Delta$ is the canonical projection;
 \item for $y\in \Delta$, $\TT^n_y$ is the common isotropy
 subgroup of all points $x\in M$ such that $J(x)=y$, \ie  $\TT^n_y=\{s\in \TT^n \tq \acta{s}{x}=x\};$
\item $\TT^n/\TT^n_y$ is the quotient
  group and $\widehat{\TT^n/\TT^n_y}$ is its Pontryagin dual. Since $\TT^n_y$ is closed in
  $\TT^n$, the dual map of $\TT^n \rOnto \TT^n/\TT^n_y$ is into:
$$\widehat{\TT^n/\TT^n_y} \rInto \widehat{\TT^n}=\ZZ^n;$$
\item For $k=(k_1,\dots,k_n)\in \ZZ^n$ and $s=(s_1,\dots,s_n)\in
  \TT^n$, we denote  the
duality bracket $$\dbracket{k}{s}=s_1^{k_1}\cdots s_n^{k_n}.$$
\item For every $k\in \ZZ^n$, we denote:
$\Delta(k)=\{y\in \Delta \tq s \in \TT^n_y \Rightarrow
\dbracket{k}{s}=1 \}.$
Then $k\in \widehat{\TT^n/\TT^n_y}$ holds
if and only if $y\in \Delta(k)$. We denote moreover:
$\Delta(\infty)=\{y\in \Delta \tq \TT^n_y=\{1\}\}$.
\end{itemize}
Finally we set $G_0$ be the bundle: 
$$G_0=\bigcup_{y\in \Delta}  \widehat{\TT^n/\TT^n_y}=\bigcup_{k\in \ZZ^n}
\Delta(k)\times \{k\} \subset
\Delta\times\ZZ^n,$$
endowed with the topology induced by $\Delta \times \ZZ^n$. 

\begin{prop}
If there exists a {\em continuous} section $\sigma$ of $J$ :
$M \pile{\lDashto^\sigma \\ \rOnto_J} \Delta,$ then $G_0\rOnto \Delta$ is an \'etale bundle of commutative groups.
\end{prop}
\begin{demo}
Since the trivial bundle of groups $\Delta\times\ZZ^n$ is clearly
\'etale, it suffices to prove that $G_0$ is open.
 So it suffices to prove that the $\Delta(k)$
are open in $\Delta$, \ie that their complements
$\Delta(k)^c$ are 
closed in $\Delta$.

First we give a characterization of elements in $\Delta(k)^c$.
For each $y\in \Delta$ define a subgroup of $\TT$: $H_y=\{\dbracket{k}{s} \tq s\in \TT^n_y\}.$
We have
$$y \in \Delta(k)^c\Leftrightarrow k\notin \widehat{\TT^n/\TT^n_y}
\Leftrightarrow H_y\neq \{1\}.$$
Moreover, recalling that Lie groups do not have small nontrivial
subgroups \cite{MZ}, in the particular case of a subgroup $H$ of $\TT$  we have:
$$(\forall t\in H, \ |t-1|<\sqrt{2}) \Leftrightarrow H=\{1\}.$$
Thus
$$y \in \Delta(k)^c\Leftrightarrow (\exists g\in\TT^n_y,\ |\dbracket{k}{s}-1|\geq \sqrt{2}).$$

Second we deduce that for each sequence $(y_i)_i$ in $\Delta(k)^c$ converging to
$y$, then  $y$ is itself in $\Delta(k)^c$. 
Since $y_i\in \Delta(k)^c$, for each $i$ there exists a $s_i\in \TT^n_{y_i}$ such
that $|\dbracket{k}{s_i}-1|\geq \sqrt{2}$.
Let $(s_{i_j})_j$ be a subsequence of $(s_i)_i$ in $\TT^n$ converging
  to a $s$. Then $|\dbracket{k}{s}-1|\geq \sqrt{2}$. Moreover, since
  $s_i\in \TT^n_{y_i}$, we have $\acta{s_{i_j}}{\sigma(y_{i_j})}=\sigma(y_{i_j})$, for
  every $j$. So, taking the limit, we get $\acta{s}{\sigma(y)}=\sigma(y)$,
  hence $s\in \TT^n_y$. Then we have proved~:
$$\exists s\in \TT^n_y,\ |\dbracket{k}{s}-1|\geq \sqrt{2}.$$
\end{demo}
As a corollary, note that the bundle map $G_0 \rTo^p \Delta$ is open.

Such continuous sections $\sigma$ may not exist; from now on we suppose that
it does and we introduce more notations which may depend implicitely of this $\sigma$:
\begin{itemize}
\item We define the map $\Delta\times \TT^n \rTo^\rho M$ by: $\rho(y,s)=\acta{s}{\sigma(y)}$. 
\item The ``dual'' bundle of $G_0$ is:
$$\displaystyle \widehat{G_0}=\bigcup_{y\in \Delta} \TT^n/\TT_y^n,$$
endowed with the biggest topology such that the canonical projection
$\Delta\times \TT^n \rOnto^{\widehat{\pi}} \widehat{G_0}$ is
continuous. 
\item Since $\rho(y,s)=\rho(y',s')$  holds if
  and only if $y=y'$  and $s=s' \mod \TT_y^n$ the map $\rho$ has a one-to-one
  quotient map $\overline{\rho}$ such that the  following diagram is commutative:
$$\begin{diagram}
\widehat{G_0} & \rTo^{\overline{\rho}} & M  \\
\uTo^{\widehat{\pi}} & \ruTo^\rho & \uTo_\alpha  \\
\Delta\times \TT^n & \rTo_{\sigma\times id_{\TT^n}} & M\times \TT^n 
\end{diagram}$$
\item For every $k\in \ZZ^n\cup\{\infty\}$, we define $M_k$ be the following subset
  of $M$:
$$M_k=\rho(\Delta(k)\times \TT^n)=J\inv
(\Delta(k))=\{\acta{s}{\sigma(y)}\tq s\in \TT^n,\ y \in \Delta(k)\}.$$
In particular, note that $M_{\infty}=J\inv(\Delta(\infty))$ is the maximal stable subset of $M$
where the action of $\TT^n$ is free.
\item  For every $k\in \ZZ^n$, then define a function $k_\sigma$ by: 
$$\begin{diagram}
 & & M_k & & & \ \ \ &  k_\sigma\left(\rho(y,s)\right)=\dbracket{k}{s}\\
\Delta(k)\times \TT^n &  \ruOnto(2,1)^\rho & \rTo^k &  \rdTo(2,1)^{k_\sigma}  & \TT  \\
\end{diagram}$$
\end{itemize}
This is the main result of the section:
\begin{theo}\label{Fourier}
For every continuous section $\sigma$ of $J$ and all $k\in
\ZZ^n$, the sets $M_k$ are open in $M$, and the functions
$k_\sigma$ are continuous. Moreover, for every continuous section $\sigma$, there is an isomorphism of
\cstar-algebras
$\Czero(M) \rTo_\sim^\isom \cstar(G_0),$
such that, for every $f\in \Cc(G_0)$ and $x\in M$:
$$(\isom\inv f)(x)=\sum_{k\in  \widehat{\TT^n/\TT^n_{y}}} f(J(x),k)
  k_\sigma(x), \ \ \ \ \hbox{ where } y=J(x).$$
\end{theo}
We decompose the proof in several lemmas:
\begin{lemme}\label{homeo}
The map $\overline{\rho}$ is an homeomorphism.
\end{lemme}
\begin{demo}
Since the group $\TT^n$ is compact, its action $\alpha$ on $M$ is a
proper map, and $\sigma$ is proper too, since, for every compact subset
$K$ of $M$, the set $\sigma\inv(K)=J(K)$ is compact. Then $\overline{\rho}\circ \widehat{\pi}
=\alpha \circ (\sigma\times id_{\TT^n})$ is a closed and continuous map. Then $\overline{\rho}$ itself
  is closed, since $\widehat{\pi}$ is continuous and
  surjective, and $\overline{\rho}$ is continuous too, by the choice of topology of
  $\widehat{G_0}$, hence $\overline{\rho}$ is a homeomorphism.
\end{demo}
Easy consequences of the lemma are:
\begin{itemize}
\item the bundle map $\widehat{G_0} \rOnto^{\widehat{p}=J\circ
    \overline{\rho}} \Delta$  is continuous and proper;
\item the bundle duality bracket $\widehat{G_0}\underset{\Delta}*G_0 \rTo \TT$
  is continuous, denoting $\widehat{G_0}\underset{\Delta}*G_0$ the fibered
  product.
\end{itemize}
Then we have:
\begin{lemme}
For every function $f\in \cstar(G_0)$ let us define its 
Fourier transform as the function:
\begin{center}\sfonction{\mathcal{F}f}{\widehat{G_0}}{\CC}{(y,s\mod
  \TT^n_y)}{{\displaystyle \sum_{k\in \widehat{\TT^n/\TT_y^n}} f(y,k)
  \dbracket{k}{s}}.}\end{center}
Then the Fourier transform is an isomorphism of \cstar-algebras $$\cstar(G_0)
  \rTo^{\mathcal{F}} \Czero(\widehat{G_0}).$$
\end{lemme}
\begin{demo}
Since $\widehat{G_0}$ is homeomorphic to $M$, the space
$\widehat{G_0}\times G_0$ is normal, and $
\widehat{G_0}\underset{\Delta}*G_0$ is  a closed subset. Then, using
Uhrysohn's Theorem, for $f\in \Cc(G_0)$ the
continuous map $(y,s,k) \rMapsto f(y,k)\dbracket{k}{s}$ defined on
$\widehat{G_0}\underset{\Delta}*G_0$ admits a continuous extension $K$
on $\widehat{G_0}\times G_0$. Then using \cite[chap X,\S3, n$^o$ 4,
th. 3]{BOUtg} and the fact that $G_0$ is \'etale hence the 
counting measure is a {\em continuous} Haar
system, we get the continuity of $\mathcal{F}f$ since:
$$(\mathcal{F}f)(y,s \mod \TT^n_y)=\sum_{k\in
  \widehat{\TT^n/\TT^n_y}} K((y,s),(y,k)).$$
Moreover, using
$$\left|\mathcal{F}f(y,s \mod \TT^n_y)\right|\leq \sum_{k\in
  \widehat{\TT^n/\TT^n_y}} \left|f(y,k)\right|,$$
  $\mathcal{F}f$ has its support contained in
the set  $\widehat{p}\inv(p(\supp f))$ which is compact, since
  $\widehat{p}$ is proper.

We consider $\Cc(G_0)$ as a dense sub-*-algebra of
$\cstar(G_0)$, and $ \Cc(\widehat{G_0})$ as a sub-*-algebra
of $\Czero(\widehat{G_0})$, dense too for the $\sup$ norm.
It is then an easy calculation to verify that the map
$\Cc(G_0)\rTo^{\mathcal{F}} \Cc(\widehat{G_0})$ is a morphism of commutative
*-algebras. 

Let us prove now, that this $\mathcal{F}$ is isometric. For $f\in \Cc(G_0)$
(resp. $\fcl \in \Cc(\widehat{G_0})$) and $y\in \Delta$, we denote
$f_y$ (resp. $\fcl_y$) its restriction to
$\widehat{\TT^n/\TT^n_y}$ (resp. $\TT^n/\TT^n_y$). Then, from the
definition of Fourier transform we have $(\mathcal{F}f)_y=
(\mathcal{F}_y)(f_y)$, where
$\cstar\left(\widehat{\TT^n/\TT^n_y}\right) \rTo^{\mathcal{F}_y}
\Continue\left(\TT^n/\TT^n_y\right)$
 is the classical isometric Fourier-Gelfand transform for commutative groups.
Since the bundle map
$G_0\rTo^p\Delta$ is open we get a continuous field of \cstar-algebras
$\left(\cstar(\widehat{\TT^n/\TT^n_y})\right)_{y\in \Delta}$, which
  \cstar-algebras of continuous sections is $\cstar(G_0)$; then, for
  every $f\in \Cc(G_0)$ we get:
$$\|f\|_{\cstar(G_0)}=\sup_{y\in \Delta}
\|f_y\|_{\cstar(\widehat{\TT^n/\TT^n_y})}=\sup_{y\in \Delta}
\|\mathcal{F}_y
(f_y)\|_{\Continue\left(\TT^n/\TT^n_y\right)}=\|\mathcal{F}f\|_{\Czero(\widehat{G_0})}.$$

Hence, by completion, $\mathcal{F}$ has an isometric extension
$\cstar(G_0)\rTo^{\mathcal{F}} \Czero(\widehat{G_0})$. Moreover this
extension is surjective, hence an isomorphism, since $\mathcal{F}(\cstar(G_0))$ is a closed
sub-*-algebra and dense, using the
Stone-Weierstrass Theorem, because, for every distinct $(y,s \mod \TT^n_y)$ and
$(y',s' \mod \TT^n_{y'})$ in $\widehat{G_0}$, it is easy to construct a $f\in
\Cc(G_0)$ such that  $(\mathcal{F}f)(y,s \mod \TT^n_y)\neq (\mathcal{F}f)(y',s' \mod \TT^n_{y'})$.
\end{demo}
Now we give the proof of the theorem:
\begin{demo}
The sets $M_k=J\inv(\Delta(k))$ are open since $J$ is continuous and the
$\Delta(k)$ are open. Moreover since, for all $s\in \TT^n$ and $y\in
\Delta(k)$:
$$k_\sigma(\rho(y,s))=\dbracket{k}{s},$$
then $k_\sigma=k\circ \overline{\rho}\inv$, hence $k_\sigma$ is continuous.

%
Then the isomorphism $\isom$ is given by:
$$\begin{diagram}
\Czero(M) & \rTo^\isom & \cstar(G_0) \\
& \rdTo(1,2)_{\overline{\rho}^*} \ldTo(1,2)^{\mathcal{F}}& &\ \ &\isom
\fcl=\mathcal{F}\inv (\fcl\circ \overline{\rho}) & \Leftrightarrow & \isom\inv f=
(\mathcal{F}f) \circ \overline{\rho} \inv \\
& \Czero(\widehat{G_0})& \end{diagram}$$
\end{demo}
\section{Deformation groupoid of a $\TT^n$-space}
Let the $\TT^n$-space $M$ is in fact a manifold, and the quotient space $\Delta$ is moreover a
locally closed subset of a Hausdorff manifold $N$ endowed with an
action $\beta$ of $\RR^n$ by diffeomorphism, 
$$\begin{diagram}
 \pile{M \  \\ \circlearrowleft _\alpha} &\ & \rTo^J & \Delta \subset
 & \pile{N \ \\ \circlearrowleft _\beta}
  \end{diagram}.$$
Then one can construct a deformation
groupoid of $M$ (for some Poisson structure) in the
following way.
For such a given action $\beta$ there is a right action of $\ZZ^n$ on
$\RR\times N$ defined, for every $\hbar\in \RR,y\in N,k\in \ZZ^n$ by
:
$$(\hbar,y).k=(\hbar,\beta_{\hbar k}(y)).$$
Let us denote $\GG=(\RR\times N)\rtimes \ZZ^n$ 
the cross-product groupoid (\cf \cite{REN}); then
$\GG$ is an \'etale Lie bundle of groupoids parametrized by $\hbar\in
\RR$. 

Using the notations $\Delta(k)$ and $\TT_y^n$ of the previous section
we set:
$$G=\{(\hbar,y,k)\in  \RR\times N \times \ZZ^n \tq y,\beta_{\hbar
  k}(y) \in \Delta(k) \hbox{ and } \TT^n_y=\TT^n_{\beta_{\hbar
  k}(y)}\}$$
Then the fiber $\{(\hbar,y,k)\in G \tq \hbar=0\}$ is exactly
  the groupoid $G_0$ defined in the previous section. 
\begin{prop}\label{etale}
Let $M$ be a Hausdorff separable manifold with an action of $\TT^n$
such that the quotient space $\Delta$ is a locally closed part of a Hausdorff
manifold $N$, which is endowed with an action $\beta$ of $\RR^n$ by diffeomorphism. Then the
previously defined $G$ is a subgroupoid of $\GG$, with space of units
$\Gzero=\RR\times \Delta$. Moreover, setting  $\G_k=\G\cap (\RR\times
N\times \{k\})$, the groupoid $G$ is \'etale
if and only if, for every $k\in \ZZ^n$, the sets
$$s(G_k)=\{(\hbar,y) \tq y,\beta_{\hbar
  k}(y) \in \Delta(k) \hbox{ and } \TT^n_y=\TT^n_{\beta_{\hbar
  k}(y)}\}$$
are open in $\RR \times \Delta$; then $G$ is a Hausdorff locally compact
   separable amenable bundle of groupoids, with continuous haar system
  and open projection $\G \rOnto \RR$.
\end{prop}

\begin{proof}
For
every groupoid $G$ with space of units $\Gzero$, recall that, for every subset $U$ of
$\Gzero$, the set $r\inv(U)\cap s\inv(U)$ is a subgroupoid of $G$, called the {\em restriction of $G$ to $U$},  denoted
$\restr{G}{U}$, and which
space of units is $U$.

For every closed subgroup $T$ of $\TT^n$, $ \widehat{\TT^n/T}$ is a
subgroup of $\ZZ^n$, hence $(\RR\times N)\rtimes \widehat{\TT^n/T}$ a
subgroupoid of $\GG$. Moreover, one can define the set
$$\Delta_T=\{y\in \Delta \tq \TT^n_y=T\}.$$
Then we get a partition: $\displaystyle\Delta=\bigcup_T \Delta_T$ and
so we can form the restriction 
$\restr{\left((\RR\times N)\rtimes \widehat{\TT^n/T}\right)}{\RR\times
  \Delta_T}$. Since, for every $k\in \ZZ^n$, we have the partition:
$$\displaystyle \Delta(k)=\bigcup_{T\ s.t.\ k\in
  \widehat{\TT^n/T}}\Delta_T,$$ then it is easy to verify that $G$ is the disjoint union 
$\displaystyle G=\bigcup_T \restr{\left( (\RR\times N)\rtimes \widehat{\TT^n/T}\right)}{\RR\times
  \Delta_T},$
hence a subgroupoid of $\GG$.

A groupoid with open source and range is \'etale if and only if it
admits a covering by open bissections. Hence a subgroupoid $G$ of an
étale groupoid $\GG$ is itself \'etale if and only if $\Gzero$ is locally closed in $\GG\zero$, and the images by the
source map $s$ (or by the range $r$) of open bissections covering $G$ are open in
$\Gzero$. For this particular case, the $G_k$ are bissections covering $G$.

The bundle map $G \rOnto^p \RR$ is open since 
the first projection $\RR\times \Delta \rOnto^{pr_1} \RR$ is open, $s$
is open, and $p=pr_1 \circ s$. Other conditions are particular cases
of standard facts for groupoids (\cf \cite{ANA}). 
\end{proof}

The main result of the section is 
\begin{theo}\label{deformationgroupoid}
Let $M$ be a Hausdorff separable manifold with a continuous action of $\TT^n$
such that the quotient space $\Delta$ is a locally closed part of a Hausdorff
manifold $N$ endowed with an action $\beta$ of $\RR^n$ by diffeomorphisms. If: 
\begin{itemize}
\item the groupoid $G$ is étale (\cf prop. \ref{etale}),
\item the projection $M \rTo^J
  N $ is smooth and has a continuous section $\sigma$, such that, for every $k\in \ZZ^n$,
  the functions $M_k \rTo^{k_\sigma}\TT$ are smooth,
\end{itemize}
then $M$ admits a Poisson
  bracket $\{.,.\}_\Omega$ such that $G$ is a deformation
  groupoid of $(M,\Omega)$.
\end{theo} 

This result is closed to the one of Rieffel \cite{Rie} up
to quite different technical conditions.
The proof consists in making use of Proposition \ref{deformation}:
we have just seen, in the previous proposition, that $G$ is a
convenient groupoid when étale; 
the condition 1 of Proposition \ref{deformation} comes from Theorem \ref{Fourier}. For the purpose
of condition 2, for every locally closed subset $Y$ of a manifold $\widetilde{Y}$,
we define the space $\Cci(Y\subset \widetilde{Y})$ to be the space of
continuous compactly supported functions $f\in \Cc(Y)$ which admit a
smooth extension $\ff\in \Ci(\widetilde{Y})$. Then we get the
following lemmas.
\begin{lemma}
With the assumptions of the Theorem \ref{deformationgroupoid}, then \\
$\mathcal{Q}=\Cci(G \subset
\GG)$ is a dense sub-*-algebra of $\cstar(G)$.
\end{lemma}

\begin{proof}
The space $\Cci(G \subset
\GG)$ is clearly a sub-vector space of $\Cc(G)$ and is
stable by involution, since, for each $f$ with extension
$\ff\in \Ci(\GG)$, the map $f^*$ has an extension
$g \rMapsto \overline{\ff(g\inv)}$ which is smooth, since the
inverse map of $\GG$ is a diffeomorphism.

Since the open bissections $G_k=G\cap (\RR\times N\times
\{k\})$, for $k\in\ZZ^n$, form a partition of $G$ we get 
$$\Cci(G\subset\GG)=\bigoplus_{k\in
  \ZZ^n} \Cci(G_k \subset \RR\times N\times
\{k\}).$$
Then, to prove the stability of $\Cci(G\subset\GG)$ by product, one
  has only to show that, for every $l,m\in \ZZ^n$ and every $f\in \Cci(G_k \subset \RR\times N\times
\{l\}),\ f'\in \Cci(G_l \subset \RR\times N\times
\{m\})$, then $f*f'$ is in $\Cci(G\subset\GG)$. This comes easily from
 the formula:
$$(f*f')(\hbar,y,k)=\begin{cases} f(\hbar, \beta_{\hbar
    l}(y),l)f'(\hbar,y,m) & \hbox{ if } k=l+m \\ 0 & \hbox{ otherwise.} \end{cases}$$

Since $\Cc(G)$ is dense in $\cstar(G)$, then $\Cci(G\subset \GG)$ is
dense in $\cstar(G)$ if every $f\in \Cc(G)$ can be approached
uniformly by
functions $f_n\in \Cci(G\subset \GG)$ since the topology induced
by $\cstar(G)$ on $\Cc(G)$ is the sup-norm topology. To prove this we
remark, since $G$ is locally closed in $\GG$, that a function
$f_n\in\Cc(G)$ is in $\Cci(G\subset\GG)$ if and only if it admits a
smooth {\em and compactly supported}
extension $\ff_n\in \Cci(\GG)$. For every $f\in \Cc(G)$, since its
support is compact, $f$ admits a continuous extension $\widetilde{f}$
on $\GG$, such that $\supp \widetilde{f}
\cap G$ is compact. Then this extension $\widetilde{f}$ can be uniformly
approximated by $\ff_n\in \Cci(\GG)$; and using smooth partition of
unity we can suppose that $\supp \ff_n \subset \supp \widetilde{f}$. Then the
restrictions of the $\ff_n$ to $G$,
$f_n=\restr{\ff_n}{G}$, 
are continuous and have their supports included
in $\supp \ff \cap G$ hence compact; so, using the previous
characterization, we have $f_n\in \Cci(G\subset \GG)$ and
:
$$\sup_{G} \left|f-f_n\right| \leq \sup_{\GG}
\left|\ff-\ff_n\right| \rTo_{n \to \infty} 0.$$
\end{proof}

As a corollary $\mathcal{Q}_0=\Cci(\G_0\subset \GG_0)$ is dense in
$\cstar(\G_0)$ (condition $2a$).  And from the formula:
$$(\isom\inv f_0)(x)=\sum_{k\in \widehat{\TT^n/\TT^n_{J(x)}}}
  f_0(J(x),k)k_\sigma(x),$$
of Theorem \ref{Fourier} and the assumptions that $J$ and the
  $k_\sigma$ are smooth, we get $\isom\inv\mathcal{Q}_0 \subset
  \Cizero(M)$ (condition $2b$). The
  condition $2c'$, hence  the Theorem \ref{deformationgroupoid} follows
  from the lemma:
\begin{lemma}\label{limitlemma}
With the assumptions of the Theorem \ref{deformationgroupoid}, for
every $f,f'\in \Cci(G\subset \GG)$ and every $\tilde{e}\in
\Cci(\GG)$ which extends $f*f'-f'*f$ on $\GG$ (resp. on a neighborhood in $\GG$
of a given point of $\G_0$
), then 
the restriction $d=\restr{\tilde{d}}{\G}$ of
$\tilde{d}\in\Ci(\tilde{G})$ defined by:
$$\widetilde{d}(\hbar,y,k)=\inverse{\imath}\int_0^1 \frac{\partial\widetilde{e}}{\partial
  \hbar}(\hbar t,y,k)\, dt,$$
\begin{enumerate}
\item is a continuous extension of the section $\hbar
  \rMapsto \inverse{\imath\hbar}[f_\hbar,{f'}_\hbar]$ on $\G$
  (resp. on a neighborhood in $\G$ of the considered point of $G_0$
).
\item has the same support as $f*f'-f'*f$ hence compact
\item has value on $\G_0$ (resp. at the considered point of $G_0$) only depending on $f_0$ and $f_0'$, and given by
 :
$$\{f_0,f_0'\}_0(0,y,k)= \inverse{\imath} \frac{\partial \tilde{e}}{\partial
  \hbar}(0,y,k).$$
\end{enumerate}
\
\end{lemma}
\begin{proof}
\begin{enumerate}
\item The 0-order Taylor integral formula gives:
$$\tilde{e}(\hbar ,y,k)=\tilde{e}(0,y,k)+\imath\hbar
\tilde{d}(\hbar ,y,k).$$
For every $(\hbar,y,k)$ in $\G$ (resp. in the considered
neighbourhood of $\G_0$), since $\cstar(\G_0)$ is commutative, and
since $(0,y,\hbar)$ is in $\G_0$, we get:
$$\tilde{e}(0,y,k)=(f*f'-f'*f)(0,y,k)=0.$$
Hence, on $\G$, we get:
$$f*f'-f'*f=\imath\hbar d,$$
\ie $d$ is a continuous extension of $\frac{f*f'-f'*f}{\imath\hbar}$.
\item Since the bundle map $\G \rOnto
  \RR$ is open, the set $G-G_0$ is dense in $G$, hence the continuous extension $d$ is unique, and the
  previous relation $f*f'-f'*f=\imath\hbar d$ shows that $f*f'-f'*f$
  and $d$ have the same support.
\item Since the map $(f,f') \rMapsto f*f'-f'*f$ is bilinear, so is
  $(f,f') \rMapsto \restr{d}{\G_0}$ too. Then it is sufficient to
  prove that $f_0=0$ implies 
  $\restr{d}{G_0}=0$. In the same way as for $\tilde{e}$ (we
  had $\restr{\tilde{e}}{G_0}=0$) we get:
$$\exists \delta f\in \Cci(G\subset\GG),\ f=\imath\hbar \delta f.$$
Then $\frac{f*f'-f'*f}{\imath\hbar}=\delta f*f'-f'*\delta f$; this
can be continously extended by 0  at $\hbar =0$.
\end{enumerate}
\end{proof}

\section{Computation of the Poisson bracket}

With some additional technical assumptions to the conditions of the Theorem
\ref{deformationgroupoid} one can compute explicitly the Poisson
structure $\Omega$ on $M$ such that $G$ is a deformation groupoid of
$(M,\Omega)$. We introduce some notations and conventions. We identify
$\RR^n$ with its Lie algebra $Lie(\RR^n)$. We fix here a basis
$E_1,\dots,E_n$ of $\RR^n$, such that we identify moreover $\RR^n$ with
$Lie(\TT^n)$, the Lie
algebra of $\TT^n$, and with its linear dual $Lie(\TT^n)^*$; in
particular $\ZZ^n=\dual{\TT^n}$ is viewed as a lattice in $\RR^n$, and
denoting $\vbracket{.}{.}$ the duality bracket of $\RR^n$, and
$Lie(\TT^n) \rTo^\exp \TT^n$ the exponential map, for every
$k=(k_1,\dots,k_n)\in \ZZ^n$ and every $X\in \RR$, we get:
$$\vbracket{E_i}{E_j}=\delta_{i,j},\ \ \vbracket{k}{E_i}=k_i,\ \
  \dbracket{k}{\exp X}=e^{\imath \vbracket{k}{X}}.$$
Moreover
  the vector field on $M$ (resp. $N$) of the
  infinitesimal action of $\alpha$ (resp. $\beta$) in the direction
  $X\in \RR^n$ will be denoted by $\xi^\alpha_X$
  (resp. $\xi^\beta_X$) \ie:
$$\forall x\in M,\ \ \xi^\alpha_X(x)=\restr{\frac{d}{d\hbar}
  \acta{\exp(\hbar X)}{x}}{\hbar=0}, \ \hbox{ and }\ \  \forall y\in N,\ \ \xi^\beta_X(y)=\restr{\frac{d}{d\hbar}
  \beta_{\hbar X}(y)}{\hbar=0}.$$
We denote $E_1^*,\dots,E_n^*$ the functions on $\G_0$ defined by:
$$E_i^*(0,y,k)=(k|E_i).$$
%

The first step is:

\begin{prop}
With the assumptions of Theorem \ref{deformationgroupoid}, let
$\Omega$ be the Poisson bivector on $M$ such that $G$ is a deformation
groupoid of $(M,\Omega)$. Then, for every
$f,f'\in \Cci(\G\subset\GG)$, every extension $\ff,\ff'\in \Cci(\GG)$
respectively of $f$ and $f'$, denoting $f_0=\restr{f}{\G_0}$ and $f_0'=\restr{f'}{\G_0}$,
and $\isom$ being the isomorphism provided by Theorem \ref{Fourier}, we
get:
{\scriptsize $$ \{\isom\inv f_0,\isom\inv f_0'\}_\Omega=\inverse{\imath} \sum_{i=1}^n
\isom\inv\left[\restr{d\ff\left(\xi^\beta_{E_i}\right)}{\G_0}\right] \left[\isom\inv(E_i^* f_0')\right]-\left[\isom\inv(E_i^*
f_0)\right] \isom\inv\left[\restr{d\ff'\left(\xi^\beta_{E_i}\right)}{\G_0}\right].$$}
\end{prop}

\begin{proof}
Since the Poisson bracket $\{.,.\}_\Omega$ on $M$ is given by the
relation:
$$\isom\inv\{f_0,f_0'\}_0=\{\isom\inv f_0,\isom\inv f_0'\}_\Omega,$$
we have only to prove that:
$$\{f_0,f_0'\}_0=\inverse{\imath} \sum_{i=1}^n
\restr{d\ff\left(\xi^\beta_{E_i}\right)}{\G_0}*(E_i^* f_0')-(E_i^*
f_0)*\restr{d\ff'\left(\xi^\beta_{E_i}\right)}{\G_0}.$$

We can restrict to functions $f,f'$ with extensions $\ff,\ff'$ both
supported in $\RR\times N\times\{l\}$ and $\RR\times N\times\{m\}$,
for some $l,m\in \ZZ^n$. Then $f*f'$ and $f'*f$ are both supported in $\RR\times N\times\{m+l\}$

Since we have supposed that $\G$ is \'etale, then every
element in $\G_0\cap (\RR\times N\times\{m+l\})$  has a neighborhood $V$ in $\G$ such that, for every
$(\hbar,y,m+l)$ in $V$, we have together:
$$(\hbar,y,l),\ (\hbar,\beta_{\hbar l}(y),m),\
(\hbar,y,m),\ (\hbar,\beta_{\hbar m}(y),l) \in G.$$
Then, on $V$, we get:
$$f*f'-f'*f=\widetilde{e},$$
where $\widetilde{e}\in \Ci(\GG)$ is the function defined by:
$$\widetilde{e}(\hbar,y,m+l)=\ff(\hbar,\beta_{\hbar
  l}(y),m)\ff'(\hbar,y,l)-\ff'(\hbar,\beta_{\hbar
  m}(y),l)\ff(\hbar,y,m).$$
Due to lemma \ref{limitlemma} we get:
$$\{f_0,f_0'\}_0(0,y,k)=\inverse{\imath} \frac{\partial \widetilde{e}}{\partial
  \hbar}(0,y,k).$$
We compute:
\begin{eqnarray*}
\frac{\partial\widetilde{e}}{\partial
  \hbar}(\hbar,y,k)&=&\left(\frac{\partial\ff}{\partial
  \hbar}(\hbar,\beta_{\hbar
  l}(y),m)+d\ff\left(\frac{d\beta_{\hbar l}}{d\hbar}(y)\right)(\hbar,\beta_{\hbar
  l}(y),m)\right)\ff'(\hbar,y,l)\\&&+\ff(\hbar,\beta_{\hbar
  l}(y),m)\frac{\partial \ff'}{\partial \hbar}(\hbar,y,l)- \ff'(\hbar,\beta_{\hbar
  m}(y),l)\frac{\partial \ff}{\partial \hbar}(\hbar,y,m)\\&&-\left(\frac{\partial\ff'}{\partial
  \hbar}(\hbar,\beta_{\hbar
  m}(y),l)+d\ff'\left(\frac{d\beta_{\hbar m}}{d\hbar}(y)\right)(\hbar,\beta_{\hbar
  m}(y),l)\right)\ff(\hbar,y,m).
\end{eqnarray*}
For $\hbar=0$, since $\beta_{0}(y)=y$ there remains
:
$$\frac{\partial \widetilde{e}}{\partial
  \hbar}(0,y,k)=d\ff(\xi^\beta_l)(0,y,m)
  \ff'(0,y,l)-d\ff'(\xi^\beta_m)(0,y,l)\ff(0,y,m).$$
Since $\displaystyle \xi_l^\beta=\sum_{i=1}^n
  \vbracket{l}{E_i} \xi_{E_i}^\beta$ we get then:
{\scriptsize $$\{f_0,f_0'\}_0(0,y,k)=\inverse{\imath}\sum_{i=1}^n
d\ff(\xi^\beta_{E_i})(0,y,m)E_i^*(0,y,l)\ff'(0,y,l)-d\ff'(\xi^\beta_{E_i})(0,y,l)E_i^*(0,y,m)\ff(0,y,m),$$}
\ie 
$$\{f_0,f_0'\}_0=\inverse{\imath} \sum_{i=1}^n
\restr{d\ff\left(\xi^\beta_{E_i}\right)}{\G_0}*(E_i^* f_0')-(E_i^*
f_0)*\restr{d\ff'\left(\xi^\beta_{E_i}\right)}{\G_0}.$$
\end{proof}

The second step is to compute $\isom\inv\restr{d\ff\left(\xi^\beta_{E_i}\right)}{\G_0}$ and $\isom\inv(E_i^*
f_0)$ with respect to $\isom\inv f_0$. We use some additional conditions
to obtain:

\begin{theo}\label{Poissonbracket}
With the assumptions of Theorem \ref{deformationgroupoid} and suppose
moreover that:
\begin{enumerate}
\item the torus $\TT^n$ acts (by $\alpha$) on $M$ by
diffeomorphisms;
\item $\Delta(\infty)=\{y \in \Delta \tq \TT_y^n=\{1\}\}$ is open in $N$
  and dense in $\Delta$ -- for density it occurs for example when the action
$\alpha$ is effective;
\item the restricted map $\Delta(\infty) \times
\TT^n \rTo^\rho M_{\infty}$ (\cf lemma \ref{homeo}) is  a diffeomorphism
\end{enumerate}
Then the Poisson structure of Theorem \ref{deformationgroupoid} is given on $M_{\infty}$ by the
bivector:
$$\Omega=\sum_{i=1}^n \xi_{E_i}^\alpha \wedge
\rho_*(\xi_{E_i}^\beta).$$
\end{theo}

\begin{proof}
%
%
Remark that, for $\mathcal{Q}=\Cci(\G\subset \GG)$, we have both:
$$\mathcal{Q}_0=\Cci(\G_0\subset \GG_0) \hbox{ and }
\isom\inv\mathcal{Q}_0 \subset \Cizero(M),$$
hence $\mathcal{Q}_0$ is dense in $\Cizero(M)$ since $\Cci(\G_0\subset
\GG_0)$ is dense in $\cstar(\G_0)$. So we have only to prove that, for every $f_0,f_0'\in
\Cci(\G_0\subset \GG_0)$:
{\footnotesize \begin{eqnarray*}
&&\{\isom\inv f_0, \isom\inv f_0'\}_\Omega=\bracket{d(\isom\inv f_0)
  \otimes d(\isom\inv f_0')}{\Omega}\\
&&= \sum_{i=1}^n d\left(\isom\inv f_0\right)\left(\xi^\alpha_{E_i}\right)d\left(\isom\inv f_0'\right)\left(\rho^*\left(\xi^\beta_{E_i}\right)\right)-d\left(\isom\inv f_0'\right)\left(\xi^\alpha_{E_i}\right)d\left(\isom\inv f_0\right)\left(\rho^*\left(\xi^\beta_{E_i}\right)\right).\end{eqnarray*}}
Moreover, it can be proved only on $M_{\infty}=J\inv(\Delta(\infty))$ since it
is dense in $M$. In fact it is a consequence of both:
$$\isom\inv(E_i^*f_0)=\inverse{\imath}d(\isom\inv
  f_0)(\xi_{E_i}^\alpha) \ \hbox{ and }\ \ \isom\inv\restr{d\ff\left(\xi^\beta_{E_i}\right)}{\G_0}=d\left(\isom\inv
  f_0\right)(\rho_*(\xi_{E_i}^\beta)).$$

For the first formula, recall from Theorem \ref{Fourier} that, for all $x\in M$: 
$$\isom\inv(f_0)(x)=\sum_{k\in \ZZ^n} f(0,J(x),k)k_{\sigma}(x).$$
Using condition (1), since, for all $x\in M$: 
$$dk_{\sigma}(\xi_{E_i}^\alpha)(x)=\restr{\frac{d}{dt}k_{\sigma}\left(\acta{\exp(tE_i)}{x}\right)}{t=0}=
\restr{\frac{d}{dt}\dbracket{k}{\exp(tE_i)}}{t=0}k_{\sigma}(x)=\imath
\vbracket{k}{E_i}k_{\sigma}(x),$$
and since $J$ is invariant with respect to the action $\alpha$, we get, for
all $x\in M$: 
\begin{eqnarray*}
d_x(\isom\inv f_0)(\xi_{E_i}^\alpha)&=&\sum_{k\in
  \ZZ^n} f(0,J(x),k)dk_{\sigma}(\xi_{E_i}^\alpha)(x)\\&=&\imath \sum_{k\in
  \ZZ^n} f(0,J(x),k)E_i^*(0,J(x),k)k_{\sigma}(x)=\imath
  \isom\inv(E_i^*f_0)(x).
\end{eqnarray*}

For the second formula, since $\xi_{E_i}^\beta$ is a vector field on
$\Delta(\infty)$, using conditions (2) and (3),  for all $s\in \TT^n$ and
$y\in \Delta(\infty)$, we get:
\begin{eqnarray*}
d_{\rho(s,y)}\left(\isom\inv
  f_0\right)(\rho_*(\xi_{E_i}^\beta))&=&d_{(s,y)}\left(\isom\inv
  f_0\circ \rho \right)(\xi_{E_i}^\beta)\\
&=& \sum_{k\in
  \ZZ^n} d_{y}\ff(\xi_{E_i}^\beta)\dbracket{k}{s}\\
&=& \left(\left(\isom\inv\restr{d\ff(\xi_{E_i}^\beta)}{\G_0}\right)\circ
  \rho\right)(s,y)
\end{eqnarray*}
hence: $d\left(\isom\inv
  f_0\right)(\rho_*(\xi_{E_i}^\beta))=\isom\inv\left(\restr{d\ff(\xi_{E_i}^\beta)}{\G_0}\right).$
\end{proof}

\section{Application to toric manifolds}

\newcommand{\group}{\texttt{G}}
\newcommand{\id}{\textrm{id}}

First we recall some facts concerning toric manifolds. We will use the
notations and conventions of sections 1 and 3.

\begin{defi}
The action $\alpha$ of a Lie group $\group$ on a symplectic manifold
$(M,\omega)$ by  diffeomorphism is {\em Hamiltonian} when there exists
a so-called {\em moment map} $M \rTo^J Lie(\group)^*$ such that:
$$\forall \fcl\in\Ci(M),\ \forall X\in Lie(\group),\
df\left(\xi^\alpha_X\right)=\{\vbracket{J}{X},\fcl\}_\omega.$$
When $M$ is connected, $J$ is unique up to a constant.

A {\em toric manifold} is a compact
and connected symplectic manifold $(M,\omega)$ endowed with an effective hamiltonian action
$\alpha$ of $\TT^n$ such that $M$ has (real) dimension $2n$.
\end{defi}

In the '80s Atiyah \cite{ATI},  Guillemin and
Sternberg \cite{GS} proved that, for a toric manifold, the image of the moment
map $J(M)$ is a convex polytope in \\$\RR^n=Lie(\TT^n)^*$. And Delzant \cite{DEL}
completed this result by the following ones:
\begin{itemize}
\item the map $J$ is $\TT^n$-invariant and the quotient map is a bijection 
  $M/\TT^n \simeq J(M)$, hence an homeomorphism. So, from now on, we identify $\Delta=M/\TT^n$
  with the polytope $J(M)$ in $N=\RR^n$.
\item Every isotropy group $\TT^n_y$ is connected and depends only on
  the open face of the polytope $\Delta$ containing $y$; more precisely there exists a parametrization of $\Delta$ such that:
$$\Delta=\{y\in \RR^n \tq \forall j\in\{1,\dots,n_0\},\
\vbracket{y}{X_j}\geq \lambda_j\}, \ \ \ X_j\in\ZZ^n,\ \ \lambda_j\in
\RR$$
$$Lie(\TT_y^n)=span\{X_j\tq \vbracket{y}{X_j}=\lambda_j\}.$$
In particular $\Delta(\infty)$ is the topological interior of $\Delta$ in $\RR^n$.
\item Two toric manifolds having the same polytope $\Delta$ are
  diffeomorphic, in a  $\TT^n$-equivariant symplectic way. So a toric
  manifold is totally characterized by its polytope.
\item When $\Delta$ is the polytope of a toric manifold $M$,
  then there is an explicit
  construction of $M$ (with action an symplectic form $\omega$) from $\Delta$.
\end{itemize}
We need to recall briefly the main steps of the latter construction; we use the previous parametrization of $\Delta$:
\begin{enumerate}
\item $\TT^{n_0}$ has a (diagonal) hamiltonian action on $\CC^{n_0}$
  endowed with its canonical symplectic structure
  $\displaystyle \omega_0=\sum_{j=1}^{n_0} dx_j\wedge dy_j$. The
  moment map is given by:
$$J_0(z_1,\dots,z_{n_0})=(\lambda_1,\dots,\lambda_{n_0})+\inverse{2}\left(|z_1|^2,\dots,|z_{n_0}|^2\right)
\in \RR^{n_0}.$$
\item There exists a group exact sequence:
$$1 \rTo \TT^d \rInto^i \TT^{n_0} \rOnto^\pi \TT^n \rOnto 1, \ \ \
(n_0=n+d)$$
such that $d\pi(F_j)=X_j$, where $F_1,\dots,F_{n_0}$ is the canonical
basis of $\RR^{n_0}$.
\item Using $Lie(\TT^p)^*=\RR^p$, let $\RR^{n_0} \rTo^{di^*}
  \RR^{d}$ be the cotangent map of $i$, and let us define $J_d=di^*\circ
  J_0$.
Then we have an inclusion map $J_d\inv\{0\} \rInto^j \CC^{n_0}$ and $J_d\inv\{0\}$ in a stable part of $\CC^{n_0}$ for the actions of
  $\TT^{n_0}$ and of $\TT^d$ -- identified with its image by $i$ in
  $\TT^{n_0}$. Moreover these two actions commute and we get a
  quotient action of $\TT^n\simeq \TT^{n_0}/\TT^d$ on
  $J_d\inv\{0\}/\TT^d$.
\item Using Marsden-Weinstein's symplectic reduction \cite{MW} the
  canonical projection $J_d\inv\{0\} \rOnto^{pr} J_d\inv\{0\}/\TT^d$ is
  a submersion, hence $J_d\inv\{0\}/\TT^d$ is a manifold of dimension $2n$, with
  symplectic structure $\omega$ such that
  $pr^*\omega=j^*\omega_0$. Moreover the action of $\TT^n$ is
  hamiltonian with moment map $J$. Calculations prove that the image
  of $J$ is then $\Delta$ so we get the following commutative
  diagram:
\begin{diagram}
 M \simeq & J_d\inv\{0\}/\TT^d& \lOnto^{pr} & J_d\inv\{0\} & \rInto^j & \CC^{n_0}\\
&& \rdOnto(1,2)_{J}    &         &
\rdTo(1,2)~{J_0\circ j}
 \ldOnto(1,2)_{J_0} &&\rdTo(1,2)^{J_d}  & & J_0\circ j=d\pi^*\circ
J \circ pr \\
&& \RR^n & \rInto_{d\pi^*} & \RR^{n_0} &
\rOnto_{di^*} & \RR^d &  & di^*\circ d\pi^*=0\\
\end{diagram} 
\end{enumerate}
The main result of this article is:
\begin{theo}\label{toric}
Let $(M,\omega)$ be a toric manifold with the action $\alpha$ of
$\TT^n$, let $M \rTo^J N=\RR^n$ be the moment map, and let $\beta$ be the
action of $\RR^n$ on $N$ by translation:
$$\beta_X(y)=y+X.$$
Then, with the same notations, the groupoid $G$ of Theorem
\ref{deformationgroupoid} is a deformation groupoid of $(M,\omega)$.
\end{theo}
We will just prove that the conditions of Theorem
\ref{deformationgroupoid} are fullfilled and then use Theorem 
\ref{Poissonbracket} to prove that the Poisson structure $\Omega$ so
obtained is the same as the one coming from $\omega$, \ie $\{.,.\}_\omega=\{.,.\}_\Omega$.

\begin{lemma} With respect to the assumptions of Theorem \ref{toric},
  for every $k\in \ZZ^n$, the following set is open in $\Delta\times
  \RR^n$:
$$s(G_k)=\{(\hbar,k)\in \RR\times \Delta \tq y,y+\hbar k\in \Delta(k) \hbox{ and }
\TT_y^n=\TT^n_{y+\hbar k}\}.$$
\end{lemma}

\begin{proof}
Recall that the isotropy subgroups $\TT^n_y$ are connected and 
$$Lie(\TT_y^n)=span\{X_j\tq \vbracket{y}{X_j}=\lambda_j\}.$$
So we get:
$$\Delta(k)=\{y\in \Delta \tq \forall j\in\{1,\dots,n_0\},\
\vbracket{y}{X_j}=\lambda_j \Rightarrow
\vbracket{k}{X_j}=0\}=\Delta\cap U(k),$$
where $U(k)$ is the open subset of $\RR^n$ defined by:
$$U(k)=\{y\in \RR^n \tq \forall j\in\{1,\dots,n_0\},\
\vbracket{k}{X_j}\neq 0\Rightarrow \vbracket{y}{X_j}>\lambda_j 
\}.$$

Since we have 
$$\vbracket{y+\hbar k}{X_j}=\vbracket{y}{X_j}+\hbar
\vbracket{k}{X_j},$$
then we get that:
$$ y\in \Delta(k) \Rightarrow Lie(\TT_y^n)\subset
Lie(\TT_{y+\hbar k}^n).$$
In the same way we get: $ y+\hbar k\in \Delta(k) \Rightarrow
Lie(\TT_{y+\hbar k}^n)\subset  Lie(\TT_y^n).$\\
Hence $y,y+\hbar k\in \Delta(k)$ implies that $\TT_y^n=
\TT^n_{y+\hbar k}$ and we have furthermore:
$$y,y+\hbar k\in \Delta(k) \Leftrightarrow y\in \delta(k) \hbox{ and } y+\hbar
k \in U(k).$$
So 
$$s(G_k)=\{(\hbar,k)\in \RR\times \Delta \tq y,\in \Delta(k) \hbox{ and }
y+\hbar k\in U(k)\}$$
is open in $\RR\times \Delta$.
\end{proof}

For the action of $\TT^{n_0}$ in $\CC^{n_0}$ we have
$J_0(\CC^{n_0})=\lambda+(\RR_+)^{n_0}=\Delta_0$, and there is a natural section $\sigma_0$
of $J_0$ given, for every $(w_1,\dots,w_{n_0})\in \lambda+(\RR_+)^{n_0}$ by: 
$$\sigma_0(w_1,\dots,w_{n_0})=\left(\sqrt{2w_1-\lambda_1},\dots,\sqrt{2w_{n_0}-\lambda_{n_0}}\right).$$

\begin{lemma}
For every toric manifold $(M,\omega)$ constructed as before, the map $\sigma=pr\circ \sigma_0\circ d\pi^*$ is well defined and is a
section of $M\rTo^J \Delta \subset \RR^n$.
\end{lemma}

\begin{proof}
Using the previous big commuting diagram and $J_0\circ \sigma_0=\id$, for every $y\in \Delta$ we get:
$$J_d((\sigma_0\circ d\pi^*)(y))=(di^*\circ d\pi^*)(y)=0,$$
hence $(\sigma_0\circ d\pi^*)(y)$ is in $J_d\inv\{0\}$, so $\sigma$ is
well defined.

The same kind of  calculus proves:
$$d\pi^* \circ (J\circ \sigma)=d\pi^*.$$
Since $\pi$ is surjective, $d\pi^*$ is into, hence $J\circ \sigma=\id.$
\end{proof}

In particular this section $\sigma$ is continuous. Moreover, since the
restriction of $\sigma_0$ to $\Delta_0(0)=\lambda+(\RR_+^*)^{n_0}$ is smooth,
we get that $\sigma$ is smooth on $\Delta(\infty)$, hence the map $\Delta(\infty)\times \TT^n
\rTo^\rho M$ is smooth.

\begin{lemma}
For every toric manifold $(M,\omega)$ constructed as before and with
respect to the previous $\sigma$, 
for all $k\in \ZZ^n$, the map $M_k
\rTo^{k_\sigma} \TT$ is smooth.
\end{lemma}

\begin{proof}
In the same way as for the action $\alpha$ of $\TT^n$ on $M$ and the
section $\sigma$, considering the diagonal action of $\TT^{n_0}$ on
$\CC^{n_0}$ and the section $\sigma_0$, one can define~: a map
$\Delta_0 \times \TT^{n_0} \rTo^{\rho_0} \CC^{n_0}$, for every
$l=(l_1,\dots,l_{n_0})\in \ZZ^{n_0}$ the sets $\Delta_0(l)$ and
$\CC_l^{n_0}=\rho_0(\Delta_0(l) \times \TT^{n_0})$, and finally the map
$\CC_l^{n_0} \rTo^{l_{\sigma_0}} \TT.$ Explicitely we get:
$$\CC_l^{n_0}=\{(z_1,\dots,z_{n_0})\in \CC^{n_0} \tq l_j\neq 0
\Rightarrow z_j\neq 0\}$$
and, since $\sigma_0(J_0(z_1,\dots,z_{n_0}))=(|z_1|,\dots,|z_{n_0}|)$:
$$l_{\sigma_0}(z_1,\dots,z_{n_0})=\prod_{j|l_j\neq 0} \left(\frac{z_j}{|z_j|}\right)^{l_j},$$
hence $l_{\sigma_0}$ is smooth on $\CC_l^{n_0}$.

For a $k\in \ZZ^n$, we apply this construction to the particular
case of $l=d\pi^*(k)\in
\ZZ^{n_0}$. 
Then one can check that the following diagram is commutative:
$$\begin{diagram}
M_k & \lOnto^{pr} & pr\inv(M_k) & \rInto^j & \CC_l^{n_0}\\
    & \rdTo(2,1)_{k_{\sigma}} & \TT & \ldTo(2,1)_{l_{\sigma_0}}&
    & \ \ \ &l_{\sigma_0}\circ j=  k_{\sigma}\circ pr
\end{diagram}$$
Since $pr$ is a submersion and $l_{\sigma_0}\circ j$ is smooth, then
$k_{\sigma}$ is smooth too.
\end{proof}

So the conditions of Theorem \ref{deformationgroupoid} are true. To
use Theorem \ref{Poissonbracket} there remains only to prove that the
map $\Delta(\infty)\times \TT^n \rTo^\rho M_{\infty}$ is a difféomorphism; we
already know that it is one-to-one and smooth. It comes from:

\begin{lemma}
For every toric manifold $(M,\omega)$ constructed as before,
and every $(y,s)\in \Delta(\infty)\times \TT^n$, identifying the tangent
bundle of $\Delta(\infty)\times \TT^n$ at $(y,s)$
with $\RR^n\times \RR^n$, we get:
$$\forall (Y,X),(Y',X')\in \RR^n\times \RR^n,\
\rho^*\omega((Y,X),(Y',X'))=\vbracket{X}{Y'}-\vbracket{Y}{X'}.$$
\end{lemma}

\begin{proof}
For every $y\in\Delta(\infty)$ and $s\in \TT^n$ we get:
$$d_{y,s}\rho(Y,X)=\xi^\alpha_X(\rho(y,s))+d_{\sigma(y)}\alpha_s\left(d_y\sigma(Y)\right).$$
So, using the relation between $\omega$ and $J$ given by 
$\omega\left(\xi^\alpha_X,T\right)=d\vbracket{J}{X}(T)$
we get:
$$\rho^*\omega((0,X),(0,X'))=\omega\left(\xi^\alpha_X,\xi^\alpha_{X'}\right)=0.$$
$$\rho^*\omega((0,X),(Y,0))=\omega\left(\xi^\alpha_X,d\alpha_s(d\sigma(Y))\right)=\vbracket{X}{Y}.$$
Then after having computed that $\sigma_0^*\omega_0=0$, using
$pr^*\omega=j^*\omega_0$ we get:
$$ \sigma^*\omega=(pr\circ \sigma_0\circ d\pi^*)^*\omega=(\sigma_0\circ
d\pi^*)^*(pr^*\omega)=(d\pi^*)^*(\sigma_0^*\omega_0)=0.$$
Hence:
$$\rho^*\omega((Y,0),(Y',0))=(\alpha_s\circ\sigma)^*\omega((Y,0),(Y',0))=\sigma^*\omega((Y,0),(Y',0))=0.$$
\end{proof}

Hence $\rho$ is a diffeomorphism since it pulls back a nondegenerated
2-form on a nondegenerated 2-form. 

Moreover, using $\xi^\beta_Y=Y$
since $\beta$ acts by translations, we get 
$$d_{\sigma(y)}\alpha_s\left(d_y\sigma(Y)\right)=\rho^*\left(\xi^\beta_Y\right)$$
Hence, it follows immediatly that the Poisson
bivector associated to $\omega$ equals the bivector $\Omega$ given in Theorem
\ref{Poissonbracket}. So Theorem \ref{toric} is proved. 

Finally let
us add a complementary result on the structure of the  \cstar-algebras
$\cstar(G_\hbar)$ which occur in the deformation:
\begin{prop}
For every $\hbar \in \RR$, every irreducible representation of
$\cstar(G_\hbar)$ has finite dimension. In particular
$\cstar(G_\hbar)$ is a \cstar-algebra of type I .
\end{prop}

\begin{demo}
For $\hbar=0$, $\cstar(\G_0)$ is commutative, hence all representations
have dimension 1. 

For $\hbar \neq 0$, we remark that the action of
$G_\hbar$ on its base $G\zero_\hbar$ has no isotropy, \ie $\G_\hbar$
is a so-called principal groupoid. Moreover every orbit of this action
is finite hence closed. But every irreducible representation of the
\cstar-algebra of a principal groupoid is supported by a closed orbit.
\end{demo}

{\bf Aknowledgments :} The author is indebted to J. Renault,
G. Skandalis and J. Bellissard for very helpful discussions and comments.

\bibliographystyle{alpha}
\bibliography{article1}

\end{document}